\documentclass[11pt, a4paper, twoside]{article}
\usepackage{amsfonts}
\usepackage{mathrsfs}
\usepackage{amsmath}
\usepackage{amssymb}
\usepackage{fancyhdr}
\usepackage{hyperref}
\numberwithin{equation}{section}

\allowdisplaybreaks

\begin{document}

\fancyhf{}

\fancyhead[EC]{Heping Liu and Manli Song}

\fancyhead[EL]{\thepage}

\fancyhead[OC]{The restriction theorem for the Grushin operators}

\fancyhead[OR]{\thepage}

\renewcommand{\headrulewidth}{0pt}
\renewcommand{\thefootnote}{\fnsymbol {footnote}}

\title{\textbf{The restriction theorem for the Grushin operators}}
{
{\footnotetext {{}\emph{Key words and phrases}: Grushin operators, scaled Hermite operators, restriction operator.}
\author{Heping Liu and Manli Song
\footnote {Corresponding author.}}
\footnotetext{The first author is supported by National Natural Science Foundation of
China under Grant \#11371036 and the Specialized Research Fund for the Doctoral Program of Higher
Education of China under Grant \#2012000110059. The second author is supported by the China Scholarship Council under Grant \#201206010098 and the Fundamental Research Funds for the Central Universities \#3102015ZY068.}
\date{}
\maketitle
\begin{abstract}
{\bf Abstract.} We study the Grushin operators acting on $\mathbb{R}^{d_1}_x \times \mathbb{R}^{d_2}_t$ and defined by the formula
\begin{equation*}
L=-\overset{d_1}{\underset{j=1}{\sum}}\partial_{x_j}^2-\left(\overset{d_1}{\underset{j=1}{\sum}}|x_j|^2\right)\overset{d_2}{\underset{k=1}{\sum}}\partial_{t_k}^2.
\end{equation*}
We establish a restriction theorem associated with the considered operators. Our result is an analogue of the restriction theorem on the Heisenberg group obtained by D. Muller.
\end{abstract}

\newtheorem{theorem}{Theorem}[section]
\newtheorem{preliminaries}{Preliminaries}[section]
\newtheorem{definition}{Difinition}[section]
\newtheorem{main result}{Main Result}[section]
\newtheorem{lemma}{Lemma}[section]
\newtheorem{proposition}{Proposition}[section]
\newtheorem{corollary}{Corollary}[section]
\newtheorem{remark}{Remark}[section]

\section[Introduction]{Introduction}
The restriction theorem for the Fourier transform plays an important role in harmonic analysis as well as in the theory of partial differential equations. The initial work on restriction theorem was given by E. M. Stein \cite{S} on $\mathbb{R}^n$. The result is stated as follows:
\begin{theorem}(Stein-Tomas) \indent Let $1 \leq p \leq \frac{2n+2}{n+3}$. Then the estimate
\begin{equation*}
||\hat{f}||_{L^2(S^{n-1})} \leq C ||f||_{L^p(\mathbb{R}^n)}
\end{equation*}
\end{theorem}
holds for all functions $f \in L^p(\mathbb{R}^n)$. \\
A simple duality argument shows that Stein-Thomas theorem is equivalent to the following:
\begin{equation}
||f*\widehat{d\sigma_r}||_{p^{'}} \leq C_r ||f||_p \label{equ:dual}
\end{equation}
holds for all $f \in \mathscr{S}(\mathbb{R}^n)$, where $d\sigma_r$ is the surface measure on the sphere with radius r.\\
Moreover, according to the Knapp example \cite{S}, the restriction theorem fails if $\frac{2n+2}{n+3}<p\leq2$.\\
From then on, the importance of the restriction theorem has become evident and various new restriction theorems has been proved. On the other hand, the restriction theorem can be generalized to many other spaces, such as Lie groups and compact manifolds (see \cite{LW}\cite{Mu}\cite{CC}\cite{SO}\cite{T1}).\\
The aim of this paper is to study the restriction theorem associated with the Grushin operators, that is,
\begin{equation*}
L=-\Delta_x-|x|^2\Delta_t,
\end{equation*}
where $(x,t)\in \mathbb{R}^{d_1}\times \mathbb{R}^{d_2}$, while $\Delta_x$, $\Delta_t$ are the corresponding partial Laplacians, and $|x|$ is the Euclidean norm of $x$. It is obvious that $L$ is self-adjoint. The operator is closely related to the scaled Hermite operators $H(a)=-\Delta_x+a^2|x|^2$. Indeed, for a Schwartz function $f$ on $\mathbb{R}^{d_1}\times \mathbb{R}^{d_2}$, let $f^\lambda(x)=\int_{\mathbb{R}^{d_2}}f(x,t)e^{i\lambda\cdot t}\,dt$ be the inverse Fourier transform of $f$ in the $t$ variable. Applying the operator $L$ to the Fourier expansion $f(x,t)=\frac{1}{2\pi}\int_{\mathbb{R}^{d_2}}f^\lambda(x)e^{-i\lambda\cdot t}\,d\lambda$, we see that
\begin{align*}
Lf(x,t)&=\frac{1}{2\pi}\int_{\mathbb{R}^{d_2}}H(|\lambda|)f^\lambda(x)e^{-i\lambda\cdot t}\,d\lambda \\
       &=\frac{1}{2\pi}\int_0^\infty a^{d_2-1}\left(\int_{S^{d_2-1}}H(a)f^{a\epsilon}(x)e^{-ia t\cdot\epsilon}\,d\sigma(\epsilon)\right)\,da
\end{align*}
Let us recall some results about the special scaled Hermite expansion. For $k \in \mathbb{N}$, the Hermite function $h_k$ of order $k$ is the function on $\mathbb{R}$ defined by
\begin{equation*}
h_k(\tau)=(2^kk!\sqrt{\pi})^{-1/2}H_k(\tau)e^{-\tau^2/2}.
\end{equation*}
Let $\nu$ be a multiindex and $x \in \mathbb{R}^{d_1}$, we define the $d_1$-dimensional Hermite functions $\Phi_\nu$ by
\begin{equation*}
\Phi_\nu(x)=\prod_{j=1}^{d_1} h_k(x_j).
\end{equation*}
The eigenfunctions of the scaled Hermite operator $H(a)$ are given by $\Phi^a_\nu=|a|^{1/4}\Phi_\nu(\sqrt{|a|}x)$ and $H(a)\Phi^a_\nu=(2k+d_1)|a|\Phi^a_\nu$. Let $P_k(a)$ stand for the projection of $L^2(\mathbb{R}^{d_1})$ onto the $k$-th eigenspace of $H(a)$. More precisely
\begin{equation*}
P_k(x)\varphi=\sum_{|\nu|=k}(\varphi, \Phi^a_\nu)\Phi^a_\nu
\end{equation*}
Then the spectral decomposition of the opearator $H(a)$ is explicitly known:
\begin{equation*}
H(a)=\underset{k=0}{\overset{\infty}{\sum}}(2k+d_1)|a|P_k(a),
\end{equation*}
Hence the spectral decomposition of the Grushin operator $L$ is given by
\begin{equation*}
Lf(x,t)=\frac{1}{2\pi}\int_0^\infty \left(\underset{k=0}{\overset{\infty}{\sum}}(2k+d_1)a^{d_2}\int_{S^{d_2-1}}P_k(a)f^{a\epsilon}(x)e^{-ia t\cdot\epsilon}\,d\sigma(\epsilon)\right)\,da
\end{equation*}
Thus the spectrum of $L$ consists of the half line $[0,\infty)$. The Grushin operator is a self-adjoint and positive operator. Furthermore, $e^{-iat\cdot \epsilon}\Phi_\nu^a(x)$ is an eigenfunction of $L$ with the eigenvalue $(2|\nu|+d_1)|a|$. Therefore, we have
\begin{align*}
Lf(x,t)&=\frac{1}{2\pi}\int_0^\infty \left[\underset{k=0}{\overset{\infty}{\sum}}\frac{\mu^{d_2}}{(2k+d_1)^{d_2}}\int_{S^{d_2-1}}P_k\left(\frac{\mu}{2k+d_1}\right)f^{\frac{\mu\epsilon}{2k+d_1}}(x)e^{-\frac{i\mu t\cdot\epsilon}{2k+d_1}}\,d\sigma(\epsilon)\right]\,d\mu\\
&=\int_0^\infty \mu \mathcal{P}_\mu f(x,t)\,d\mu
\end{align*}
where $\mathcal{P}_\mu f(x,t)=\frac{1}{2\pi} \underset{k=0}{\overset{\infty}{\sum}}\frac{\mu^{d_2-1}}{(2k+d_1)^{d_2}}\int_{S^{d_2-1}}P_k\left(\frac{\mu}{2k+d_1}\right)
f^{\frac{\mu\epsilon}{2k+d_1}}(x)e^{-\frac{i\mu t\cdot\epsilon}{2k+d_1}}\,d\sigma(\epsilon)$
is an eigenfunction of $L$ with the eigenvalue $\mu$. $\mathcal{P}_\mu$ is called the restriction operator.\\
Let $L_k^\delta$ be the Laguerre polynomial of type $\delta$ and degree $k$ defined by
\begin{equation*}
L_k^\delta(\tau)=\frac{1}{k!}\frac{d^k}{d\tau^k}(e^{-\tau}x^{k+\delta})e^\tau x^{-\delta}, \quad \forall k \in \mathbb{N}, \delta>-1
\end{equation*}
We define the normalized Laguerre functions by
\begin{equation*}
\mathcal{L}_k^\delta(\tau)=\left(\frac{\Gamma(k+1)}{\Gamma(k+\delta+1)}\right)^{\frac{1}{2}}e^{-\frac{\tau}{2}}\tau^{\frac{\delta}{2}}L_k^\delta(\tau)
\end{equation*}
and set the Laguerre functions $\varphi_k(z)=L_k^{d_1-1}(\frac{1}{2}|z|^2)e^{-\frac{1}{4}|z|^2}, \forall z \in \mathbb{C}^{d_1}$.\\
Next we introduce the Weyl transform. Let $H^{d_1}$ be the $(2d_1+1)$-dimensional Heisenberg group. For each $a \in \mathbb{R}\backslash\{0\}$, $(z,s) \in H^{d_1}$, there is an infinite dimensional representation $\pi_a(z,s)$ which in the Schr$\ddot{o}$dinger realization acts on $L^2(\mathbb{R}^{d_1})$ in the following way. For each $\varphi \in L^2(\mathbb{R}^{d_1})$, $z=x+iy$
\begin{equation*}
\pi_a(z,s)\varphi(\xi)=e^{ias}e^{ia(x\cdot\xi+\frac{1}{2}x\cdot y)}\varphi(\xi+y).
\end{equation*}
For any integrable function $g$ on $\mathbb{C}^{d_1}$, we define the Weyl transform of $g$ by
\begin{equation*}
W_a(g)=\int_{\mathbb{C}^{d_1}}g(z)\pi_a(z,0)\,dz
\end{equation*}
For each $\varphi,\psi \in L^2(\mathbb{R}^{d_1})$,
\begin{equation*}
|(W_a(g)\varphi, \psi)|=|\int_{\mathbb{C}^{d_1}}g(z)(\pi_a(z,0)\varphi, \psi)\,dz|\leq||g||_1||\varphi||_2||\psi||_2
\end{equation*}
This shows that $W_a(g)$ is a bouded operator on $L^2(\mathbb{R}^{d_1})$ with $||W_a(g)||\leq ||g||_1$. Furthermore, from the explicit description of $\pi(z,s)$ we see that
\begin{equation*}
W_a(g)\varphi(\xi)=\int_{\mathbb{C}^{d_1}}e^{ia(x\cdot\xi+\frac{1}{2}x\cdot y)}g(x,y)\varphi(\xi+y)\,dxdy
\end{equation*}
From the above it follows that $W_a(g)$ is an integral operator whose kernel $K_g^a(x,y)$ is given by
\begin{equation*}
K_g^a(x,y)=\int_{\mathbb{R}^{d_1}}g(\xi,y-x)e^{i\frac{a}{2}\xi\cdot(x+y)}\,d\xi
\end{equation*}
\\Our main result is the following theorem:
\begin{theorem}\indent \label{restriction}For $1\leq p\leq\frac{2(d_2+1)}{d_2+3}$ and $1\leq q\leq 2\leq r\leq\infty$, the following inequality holds
\begin{equation*}
||\mathcal{P}_\mu f||_{L_x^rL_t^{p^{'}}}\leq C\mu^{2d_2(\frac{1}{p}-\frac{1}{2})+\frac{d_1}{2}(\frac{1}{q}-\frac{1}{r})-1}||f||_{L_x^qL_t^p}
\end{equation*}
for any Schwartz function $f$ on $\mathbb{R}^{d_1}\times\mathbb{R}^{d_2}$ and $\mu>0$.\\
\end{theorem}
The theorem is stated in terms of the mixed Lebesgue norm
\begin{equation*}
||f||_{L_x^qL_t^p}=\left(\int_{\mathbb{R}^{d_2}}\left(\int_{\mathbb{R}^{d_1}}|f(x,t)|^q\,dx\right)^{\frac{p}{q}}\,dt\right)^\frac{1}{p}, 1\leq p,q\leq\infty
\end{equation*}
(with the obvious modifications when $p$ or $q$ are equal to $\infty$).
\section[Restriction theorem]{Restriction theorem}
To prove the above theorem, first we state some asymptotic properties of the normalized Laguerre functions $\mathcal{L}_k^\delta(r)$ (see \cite{T2}). We let $\nu=4k+2\delta+2$ and assume $\delta>-1$.
\begin{lemma}\indent The Laguerre functions satisfy
\begin{equation*}
|\mathcal{L}_k^\delta(\tau)|\leq C\left\{
\begin{array}{ll}
(\tau\nu)^{\delta/2},& 0\leq \tau\leq1/\nu \\
(\tau\nu)^{-1/4},& 1/\nu\leq \tau\leq\nu/2 \\
\nu^{-1/4}(\nu^{1/3}+|\nu-\tau|)^{-1/4},& \nu/2\leq \tau\leq 3\nu/2 \\
e^{-\gamma \tau}, &\tau\geq 3\nu/2\\
\end{array} \right.
\end{equation*}
where $\gamma>0$ is a fixed constant.
\end{lemma}
Using the above estimates for the Laguerre functions, we can get a lower bound for the $L^1$ norm of them.
\begin{lemma}\indent \label{Laguerre}
\begin{equation*}
 \int_0^\infty \left|\mathcal{L}_k^{d_1-1}(\tau)\right|\tau^{-\frac{1}{2}}\,d\tau\leq C
\end{equation*}
\end{lemma}
Moreover, there is an interesting result which connects the Laguerre function $\varphi_{k,a}^{d_1-1}(z)=\varphi_k^{d_1-1}(\sqrt{a}z)$ with the spectral projection $P_k(a)$ (see \cite{T3}).
\begin{lemma}\indent
\begin{equation*}
W_a(\varphi_{k,a}^{d_1-1})=(2\pi)^{d_1}|a|^{-d_1}P_k(a)
\end{equation*}
\end{lemma}
In order to prove the restriction theorem, we need the estimates of the projections $\varphi\rightarrow P_k(a)\varphi$ which are given in the following proposition.
\begin{proposition}\indent \label{Hermite}For $\varphi \in L^q(\mathbb{R}^{d_1})$, $1\leq q\leq2$,
\begin{equation*}
||P_k(a)\varphi||_2\leq C|a|^{\frac{d_1}{2}(\frac{1}{q}-\frac{1}{2})}(2k+d_1)^{\frac{d_1-1}{2}(\frac{1}{q}-\frac{1}{2})}||\varphi||_q
\end{equation*}
\end{proposition}
{\bf Proof.} As $||P_k(a)\varphi||_2\leq||\varphi||_2$, it is enough to prove the above estimate when $q=1$. Since
\begin{equation*}
||P_k(a)\varphi||_2^2=(P_k(a)\varphi, P_k(a)\varphi)=(P_k(a)\varphi, \varphi)\leq||P_k(a)\varphi||_{q^{'}}||\varphi||_q,
\end{equation*}
it is enough to show that
\begin{equation}\label{infty-1}
||P_k(a)\varphi||_\infty\leq|a|^{\frac{d_1}{2}}(2k+d_1)^{\frac{d_1-1}{2}}||\varphi||_1.
\end{equation}
To prove \eqref{infty-1} we use the fact that $P_k(a)=(2\pi)^{-d_1}|a|^{d_1}W_a(\varphi_{k,a}^{d_1-1})$. This shows that $P_k(a)$ is an integral operator with the kernel $F_{k,a}(x,y)$ given by
\begin{align*}
F_{k,a}(x,y)&=(2\pi)^{-d_1}|a|^{d_1}\int_{\mathbb{R}^{d_1}}e^{i\frac{a}{2}\xi\cdot(x+y)}\varphi_{k,a}^{d_1-1}(\xi,x-y)\,d\xi\\
            &=(2\pi)^{-d_1}|a|^{d_1}\int_{\mathbb{R}^{d_1}}e^{i\frac{a}{2}\xi\cdot(x+y)}L_k^{d_1-1}\left(\frac{|a|}{2}(|\xi|^2+|x-y|^2)\right)e^{-\frac{|a|}{4}(|\xi|^2+|x-y|^2)}\,d\xi
\end{align*}
Therefore, we have the estimate
\begin{align*}
\left|F_{k,a}(x,y)\right|&\leq(2\pi)^{-d_1}|a|^{d_1}\int_{\mathbb{R}^{d_1}}\left|L_k^{d_1-1}\left(\frac{|a|}{2}(|\xi|^2+|x-y|^2)\right)\right|e^{-\frac{|a|}{4}(|\xi|^2+|x-y|^2)}\,d\xi\\
&\leq C|a|^{d_1}\int_0^\infty\left|L_k^{d_1-1}\left(\frac{|a|}{2}(r^2+|x-y|^2)\right)\right|e^{-\frac{|a|}{4}(r^2+|x-y|^2)}r^{d_1-1}\,dr\\
&\leq C|a|^{\frac{d_1}{2}}\int_0^\infty\left|L_k^{d_1-1}\left(\frac{1}{2}r^2\right)\right|e^{-\frac{1}{4}r^2}r^{d_1-1}\,dr\\
&\leq C|a|^{\frac{d_1}{2}}(2k+d_1)^{\frac{d_1-1}{2}}\int_0^\infty\left|\mathcal{L}_k^{d_1-1}(r)\right|r^{-\frac{1}{2}}\,dr\\
\end{align*}
Using the estimate of Lemma \ref{Laguerre}, we get
\begin{equation*}
\left|F_{k,a}(x,y)\right|\leq C|a|^{\frac{d_1}{2}}(2k+d_1)^{\frac{d_1-1}{2}}.
\end{equation*}
This proves \eqref{infty-1} and hence the proposition.\\\\
Now we give a proof of Theorem \ref{restriction}.\\
{\bf Proof of Theorem \ref{restriction}}.  In order to simplify the notations, we write $f$ as it were the product of two functions, that is $f(x,t)=g(x)h(t)$, with $f$ and $g$ Schwartz functions. However, in the proof we will never use this fact. We take $\alpha: \mathbb{R}^{d_1}\rightarrow \mathbb{C}$ and $\beta: \mathbb{R}^{d_2}\rightarrow \mathbb{C}$, $\alpha \in \mathscr{S}(\mathbb{R}^{d_1})$,$\beta \in \mathscr{S}(\mathbb{R}^{d_2})$. Because the spectral projections associated to the scaled Hermite operator are orthogonal, we have
\begin{align*}
\langle\mathcal{P}_\mu f, \alpha \otimes \beta \rangle &=\int_{\mathbb{R}^{d_1}}\int_{\mathbb{R}^{d_2}}\mathcal{P}_\mu f(x,t)\overline{\alpha(x) \beta(t)}\,dxdt\\
&=\frac{1}{2\pi}\int_{\mathbb{R}^{d_1}}\int_{\mathbb{R}^{d_2}}\left[\underset{k=0}{\overset{\infty}{\sum}}\frac{\mu^{d_2-1}}{(2k+d_1)^{d_2}}\int_{S^{d_2-1}}\hat{h}\left(\frac{\mu\epsilon}{2k+d_1}\right)
\right.\\
&\left.\quad \quad \quad \quad \quad \times P_k\left(\frac{\mu}{2k+d_1}\right)g(x)e^{-\frac{i\mu t\cdot\epsilon}{2k+d_1}}\,d\sigma(\epsilon)\right]\overline{\alpha(x) \beta(t)}\,dxdt\\
&=\frac{1}{2\pi}\underset{k=0}{\overset{\infty}{\sum}}\frac{\mu^{d_2-1}}{(2k+d_1)^{d_2}}\int_{\mathbb{R}^{d_1}}\int_{S^{d_2-1}}\hat{h}\left(\frac{\mu\epsilon}{2k+d_1}\right)\\
&\quad \quad \quad \quad \quad \times P_k\left(\frac{\mu}{2k+d_1}\right)g(x)
\left(\int_{\mathbb{R}^{d_2}}\overline{\alpha(x) \beta(t)}e^{-\frac{i\mu t\cdot\epsilon}{2k+d_1}}\,dt\right)\,dxd\sigma(\epsilon)\\
&=\frac{1}{2\pi}\underset{k=0}{\overset{\infty}{\sum}}\frac{\mu^{d_2-1}}{(2k+d_1)^{d_2}}\int_{S^{d_2-1}}\hat{h}\left(\frac{\mu\epsilon}{2k+d_1}\right)\overline{\hat{\beta}\left(\frac{\mu\epsilon}{2k+d_1}\right)}\,d\sigma(\epsilon)\langle P_k\left(\frac{\mu}{2k+d_1}\right)g, \alpha\rangle\\
&=\frac{1}{2\pi}\underset{k=0}{\overset{\infty}{\sum}}\frac{\mu^{d_2-1}}{(2k+d_1)^{d_2}}\int_{S^{d_2-1}}\hat{h}\left(\frac{\mu\epsilon}{2k+d_1}\right)\overline{\hat{\beta}\left(\frac{\mu\epsilon}{2k+d_1}\right)}\,d\sigma(\epsilon)\\
&\quad \quad \quad \quad \quad \times \langle P_k\left(\frac{\mu}{2k+d_1}\right)g, P_k\left(\frac{\mu}{2k+d_1}\right)\alpha\rangle\\
\end{align*}
Applying the H$\ddot{o}$lder's inequality to the inner integral we deduce that
\begin{align*}
\langle\mathcal{P}_\mu f, \alpha \otimes \beta \rangle
&\leq\frac{1}{2\pi}\underset{k=0}{\overset{\infty}{\sum}}\frac{\mu^{d_2-1}}{(2k+d_1)^{d_2}}\left(\int_{S^{d_2-1}}\left|\hat{h}\left(\frac{\mu\epsilon}{2k+d_1}\right)\right|^2\,d\sigma(\epsilon) \right)^{\frac{1}{2}} \left(\int_{S^{d_2-1}}\left|\hat{\beta}\left(\frac{\mu\epsilon}{2k+d_1}\right)\right|^2\,d\sigma(\epsilon) \right)^{\frac{1}{2}}\\
&\quad \quad \quad \times \left|\left|P_k\left(\frac{\mu}{2k+d_1}\right)
g\right|\right|_{L_x^2}\left|\left|P_k\left(\frac{\mu}{2k+d_1}\right)
\alpha\right|\right|_{L_x^2}\\
\end{align*}
By Proposition \ref{Hermite}, we have for any $1\leq q\leq 2\leq r\leq\infty$
\begin{align}
\left|\left|P_k\left(\frac{\mu}{2k+d_1}\right)g\right|\right|_{L_x^2}&\leq C\mu^{\frac{d_1}{2}(\frac{1}{q}-\frac{1}{2})}(2k+d_1)^{-\frac{1}{2}(\frac{1}{q}-\frac{1}{2})}||g||_{L_x^q} \label{g}\\
\left|\left|P_k\left(\frac{\mu}{2k+d_1}\right)\alpha\right|\right|_{L_x^2}&\leq C\mu^{\frac{d_1}{2}(\frac{1}{r^{'}}-\frac{1}{2})}(2k+d_1)^{-\frac{1}{2}(\frac{1}{r^{'}}-\frac{1}{2})}||\alpha||_{L_x^{r^{'}}}\label{alpha}
\end{align}
For $1\leq p\leq\frac{2(d_2+1)}{d_2+3}$, it follows from the restriction theorem of the Fourier transform on $S^{d_2-1}$ that
\begin{align}
\left(\int_{S^{d_2-1}}\left|\hat{h}\left(\frac{\mu\epsilon}{2k+d_1}\right)\right|^2\,d\sigma(\epsilon) \right)^{\frac{1}{2}} &\leq C\left(\frac{2k+d_1}{\mu}\right)^{d_2(1-\frac{1}{p})}||h||_{L_t^p}  \label{h}\\
\left(\int_{S^{d_2-1}}\left|\hat{\beta}\left(\frac{\mu\epsilon}{2k+d_1}\right)\right|^2\,d\sigma(\epsilon) \right)^{\frac{1}{2}} &\leq C\left(\frac{2k+d_1}{\mu}\right)^{d_2(1-\frac{1}{p})}||\beta||_{L_t^p} \label{beta}
\end{align}
Therefore, by \eqref{g}, \eqref{alpha}, \eqref{h} and \eqref{beta} we have
\begin{equation*}
\langle\mathcal{P}_\mu f, \alpha\otimes\beta \rangle\leq C\underset{k=0}{\overset{\infty}{\sum}}(2k+d_1)^{-2d_2(\frac{1}{p}-\frac{1}{2})-\frac{1}{2}(\frac{1}{q}-\frac{1}{r})}\mu^{2d_2(\frac{1}{p}-\frac{1}{2})+\frac{d_1}{2}(\frac{1}{q}-\frac{1}{r})-1}||f||_{L_x^qL_t^p}||\alpha\otimes\beta||_{L_x^{r^{'}}L_t^p}
\end{equation*}
If $d_2\geq 2$ or $1\leq q<2\leq r\leq \infty$ or $1\leq q\leq2<r\leq\infty$, because of $1\leq p\leq \frac{2(d_2+1)}{d_2+3}$, we have $2d_2(\frac{1}{p}-\frac{1}{2})+(\frac{1}{q}-\frac{1}{2})>1$. Hence, the above sum converges and consequently we have
\begin{equation*}
||\mathcal{P}_\mu f||_{L_x^rL_t^{p^{'}}}\leq C\mu^{2d_2(\frac{1}{p}-\frac{1}{2})+\frac{d_1}{2}(\frac{1}{q}-\frac{1}{r})-1}||f||_{L_x^qL_t^p}
\end{equation*}
If $ d_2=1$ and $r=q=2$, we have
\begin{equation*}
\mathcal{P}_\mu f(x,t)=\frac{1}{2\pi}\underset{k=0}{\overset{\infty}{\sum}}\frac{1}{2k+d_1}\left[P_k\left(\frac{\mu}{2k+d_1}\right)
f^{\frac{\mu}{2k+d_1}}(x)e^{-\frac{i\mu t}{2k+d_1}}+P_k\left(\frac{\mu}{2k+d_1}\right)
f^{-\frac{\mu}{2k+d_1}}(x)e^{\frac{i\mu t}{2k+d_1}}\right].
\end{equation*}
Since the operators $P_k(a)$ are orthogonal projections, we have
\begin{align*}
||\mathcal{P}_\mu f||_{L_x^2L_t^\infty}&\leq\frac{1}{2\pi}\underset{k=0}{\overset{\infty}{\sum}}\frac{1}{2k+d_1}\left[||P_k\left(\frac{\mu}
{2k+d_1}\right) f^{\frac{\mu}{2k+d_1}}||_{L_x^2}+||P_k\left(\frac{\mu}{2k+d_1}\right) f^{\frac{\mu}{2k+d_1}}||_{L_x^2}\right]\\
&\leq\frac{1}{2\pi}\underset{k=0}{\overset{\infty}{\sum}}\left[||P_k\left(\frac{\mu}{2k+d_1}\right)
f^{\frac{\mu}{2k+d_1}}||_{L_x^2}+||P_k\left(\frac{\mu}{2k+d_1}\right)
f^{\frac{\mu}{2k+d_1}}||_{L_x^2}\right]\\
&\leq\frac{1}{2\pi}\left(||f^{\frac{\mu}{2k+d_1}}||_{L_x^2}+||f^{-\frac{\mu}{2k+d_1}}||_{L_x^2}\right)\\
&\leq \frac{1}{\pi}||f||_{L_x^2L_t^1}
\end{align*}
\section[Sharpness of the range p]{Sharpness of the range p}
In this section we only give an example to show that the range of $p$ in the restriction theorem associated with sublaplacian is sharp. The example is constructed similarly to the counterexample of M$\ddot{u}$ller \cite{Mu}, which shows that the estimates between Lebesgue spaces for the operators $\mathcal{P}_\mu$ are necessarily trivial.\\
\indent Let $\varphi \in C_c^{\infty}(\mathbb{R}^m)$ be a radial function, such that $\varphi(a)=\psi(|a|)$, where $\psi \in C_c^{\infty}(\mathbb{R})$, $\psi=1$ on a neighborhood of the point $\frac{1}{d_1}$ and $\psi=0$ near 0. Let $h$ be a Schwartz function on $\mathbb{R}^{d_2}$ and define
\begin{equation*}
f(x,t)=\int_{\mathbb{R}^{d_2}} \varphi(\lambda)\hat{h}(\lambda) e^{-\frac{|\lambda|}{2}|x|^2} e^{-i\langle \lambda,t\rangle} |\lambda|^n \, d\lambda
\end{equation*}
Denote $g(x,t)=\int_{\mathbb{R}^{d_2}} \varphi(\lambda) e^{-\frac{|\lambda|}{2}|x|^2} e^{-i\langle \lambda,t\rangle} |\lambda|^n \, d\lambda
           =\int_{\mathbb{R}^{d_1+d_2}}\varphi(\lambda) e^{-\frac{|\xi|^2}{2|\lambda|}} e^{-i(\langle \lambda,t \rangle+\langle \xi,x\rangle)} \, d\xi d\lambda$\\
Hence $\widehat{g(\xi,a)}=\varphi(a) e^{-\frac{|\xi|^2}{2|a|}}$, which shows that $\hat{g}$ and consequently $g$ are Schwartz functions. On the other hand, we have $f=h*_tg$, where $"*_t"$ denotes the involution about the second variable. Then,\\
\begin{equation*}
\begin{split}
f(x,t)&=\int_0^{+\infty}\biggl(\lambda^{d_1+d_2-1}\psi(\lambda) e^{-\frac{\lambda}{2}|x|^2}\underset{S^{m-1}}{\int} \hat{h}(\lambda w)e^{-i \lambda \langle w,t
      \rangle} \, d\sigma(w) \biggr)d\lambda \\
      &=\int_0^{+\infty} \biggl( {d_1}^{-d_1-d_2}\mu^{d_1+d_2-1} \psi(\frac{\mu}{d_1}) e^{-\frac{\mu |x|^2}{2d_1}} \underset{S^{d_2-1}}{\int} \hat{h}(\frac{\mu w}{d_1}) e^{-i\frac{\mu}{d_1}\langle w,t \rangle} \, d\sigma(w) \biggr) d\mu\\
      &=\int_0^{+\infty} \mathcal{P}_\mu f(x,t) \, d\mu\\
\end{split}
\end{equation*}
where
\begin{equation*}
\mathcal{P}_\mu f(x,t)= {d_1}^{-d_1-d_2}\mu^{d_1+d_2-1} \psi(\frac{\mu}{d_1}) e^{-\frac{\mu |x|^2}{2d_1}} \underset{S^{d_2-1}}{\int} \hat{h}(\frac{\mu w}{d_1}) e^{-i\frac{\mu}{d_1}\langle w,t \rangle} \, d\sigma(w)
\end{equation*}
and it satisfies $L(\mathcal{P}_\mu f)=\mu \mathcal{P}_\mu f$.\\
Therefore, specially let $\mu=1$, we have
\begin{equation*}
\begin{split}
\mathcal{P}_1f(x,t)&={d_1}^{-d_1-d_2} e^{-\frac{|x|^2}{2d_1}} \underset{S^{d_2-1}}{\int} \hat{h}(\frac{w}{d_1}) e^{-i\frac{\langle w,t \rangle}{d_1}} \, d\sigma(w)\\
                   &={d_1}^{-d_1-d_2} e^{-\frac{|x|^2}{2d_1}} h*\widehat{d\sigma_{\frac{1}{d_1}}}(t)\\
\end{split}
\end{equation*}
\indent From the restriction theorem associated the Grushin operators, we have the estimate $||\mathcal{P}_1f||_{L_x^{r}L_t^{p^{'}}}\leq C ||f||_{L_x^qL_t^p}$.\\
Because of
\begin{equation}
||\mathcal{P}_1f||_{L_x^{r}L_t^{p^{'}}}=C||h*\widehat{d\sigma_{\frac{1}{d_1}}}||_{L_t^{p^{'}}}
\end{equation}
 and
\begin{equation}
||f||_{L_x^qL_t^p}\leq ||h||_{L_t^p}||g||_{L_x^qL_t^1}\lesssim||h||_{L_t^p}
\end{equation}
we have $||h*\widehat{d\sigma_{\frac{1}{d_1}}}||_{L_t^{p^{'}}}\leq C||h||_{L_t^p}$.\\
From the sharpness of Stein-Tomas theorem which is guaranteed by the Knapp counterexample, it would imply $p\leq \frac{2(d_2+1)}{d_2+3}$. Hence the range of $p$ can not be extended.\\

{\bf Acknowledgements:}
The work is performed while the second author studies as a joint Ph.D. student in the mathematics department of Christian-Albrechts-Universit$\ddot{a}$t zu Kiel. She wishes to express her thanks to Professor Detlef M$\ddot{u}$ller for his assistance and generous discussions on restriction theorems.

\begin{flushleft}
\vspace{0.3cm}\textsc{Heping Liu\\School of Mathematical
Sciences\\Peking University\\Beijing 100871\\People's
Republic of China\\}
\emph{E-mail address}: \text{hpliu@pku.edu.cn}
\vspace{0.3cm}\textsc{Manli Song\\School of Natural and Applied
Sciences\\Northwestern Polytechnical University\\Xi'an, Shaanxi 710129\\People's
Republic of China\\}
\emph{E-mail address}: \text{mlsong@nwpu.edu.cn}
\end{flushleft}

\end{document}